\documentclass[preprint,12pt]{elsarticle}

\usepackage{amssymb}

\begin{document}
\def\E{\ensuremath{{\bf E}}}
\def\R{\ensuremath{{\bf R}}}
\def\Z{\ensuremath{{\bf Z}}}
\def\Q{\ensuremath{{\bf Q}}}
\def\N{\ensuremath{{\bf N}}}
\def\Q{\ensuremath{{\bf D}}}
\def\P{\ensuremath{{\bf P}}}
\def\lh{\ensuremath{\underline{h}}}
\def\uh{\ensuremath{\overline{h}}}
\def\lp{\left(}
\def\rp{\right)}
\def\dist{\mbox{\rm dist}}
\def\dim{\mbox{\rm dim}}
\def\hau{\ensuremath{\mathcal{H}}}
\def\supp{\ensuremath{\mbox{\rm supp}}}
\def\rmd{\mathrm{d}}
\def\rme{\mathrm{e}}
\def\rmi{\mathrm{i}}
\newcommand{\floor}[1]{{\ensuremath{[#1]}}}
\newtheorem{thm}{Theorem}
\newtheorem{pro}[thm]{Proposition}
\newtheorem{lem}[thm]{Lemma}
\newtheorem{cor}[thm]{Corollary}
\newdefinition{defi}{Definition}
\newdefinition{rem}{Remark}
\newproof{pf}{Proof}

\def\BRem{\begin{rem}}
\def\ERem{\end{rem}}
\def\BThe{\begin{thm}}
\def\EThe{\end{thm}}
\def\BDef{\begin{defi}}
\def\EDef{\end{defi}}
\def\BPro{\begin{pro}}
\def\EPro{\end{pro}}
\def\BLem{\begin{lem}}
\def\ELem{\end{lem}}
\def\BCor{\begin{cor}}
\def\ECor{\end{cor}}
\def\BProof{\begin{pf}}
\def\EProof{\end{pf}}

\begin{frontmatter}

\title{A wavelet characterization for the upper global H\"older index}
\author[mc]{M.~Clausel}
\ead{marianne.clausel@insa-lyon.fr}
\address[mc]{Universit\'e de Lyon, CNRS, INSA de Lyon, Institut Camille Jordan UMR 5208, B\^atiment L. de Vinci, 20 av. Albert Einstein, F-69621 Villeurbanne Cedex, France.}
\author[sn,ca]{S.~Nicolay\corref{cor}}
\ead{S.Nicolay@ulg.ac.be}
\address[sn]{Universit\'e de Li\`ege, Institut de Math\'ematique, Grande Traverse, 12, B\^atiment B37, B-4000 Li\`ege (Sart-Tilman), Belgium.}
\cortext[ca]{Corresponding author. Phone:~+32(0)43669433. Fax:~+32(0)43669547.}

\begin{abstract}
In this paper, we give a wavelet characterization of the upper global H\"older index, which can be seen as the irregular counterpart of the usual global H\"older index, for which a wavelet characterization is well-known.
\end{abstract}


\end{frontmatter}

\section{Introduction}
One of the most popular concept of uniform regularity is the uniform H\"older regularity, defined from the uniform H\"older spaces $C^\alpha(\R^d)$. For any $\alpha\in (0,1)$, a bounded function $f$ belongs to $C^\alpha(\R^d)$ if there exists $C,R>0$ such that
\[
 \sup_{|x-y|\le r} |f(x)-f(y)| \le C r^\alpha
\]
for any $r\in [0,R]$. This notion can be generalized for exponents greater than one (see section~\ref{sec:lwhi}). It has been widely used to study smoothness properties of classical models such as trigonometric series (see e.g.~\cite{Wei59,Kah86}) and sample paths properties of processes (amongst these processes, let us cite the Brownian motion (see \cite{Khin24}) and the fractional Brownian motion).

In many classical cases, the smoothness behavior of the investigated model is very simple. The studied function $f$ is both uniformly H\"older and uniformly anti-H\"older (see~\cite{CN10a} and \cite{CN10b} for more details) and its smoothness properties can be characterized using a single index,
\[
\mathcal{H}=\lim_{r\to 0}\frac{\log \sup_{|x-y|\le r}|f(x)-f(y)|}{\log r}.
\]
There are many well-known examples of such models (see~\cite{Wei59,Kah86,BH99,BH00,Heurt03,Heurt05} for trigonometric series and \cite{Ber72,Ber73,Adl81,Xiao97,Xiao05} for sample paths of the FBM or some of its extensions).

Nevertheless, the smoothness properties of the model can be much more complex: in many cases, the uniform modulus of smoothness $\omega_f^1$ of $f$, that is
the map
\[
 \omega_f^1: r\mapsto \sup_{|x-y|\le r} |f(x)-f(y)|,
\]
is quite general. This is for example the case with the $\phi$--SNLD Gaussian models (see~\cite{Xiao97,Xiao05}) or the lacunary fractional Brownian motion (see~\cite{Clau10}), for which the uniform modulus of smoothness may be a general function that is not possible to estimate. It is then more convenient to describe the smoothness properties of the model using two indices:
\begin{equation}\label{eq:hlow}
\underline{\mathcal{H}}=\liminf_{r\to 0} \frac{\log \sup_{|x-y|\le r}|f(x)-f(y)|}{\log r}
\end{equation}
and
\begin{equation}\label{eq:hup}
 \overline{\mathcal{H}}=\limsup_{r\to 0}\frac{\log \sup_{|x-y|\le r}|f(x)-f(y)|}{\log r},
\end{equation}
related to the behavior of the uniform modulus of smoothness of $f$ near $0$.

Even in the case of Gaussian models, the estimation of these two indices is still an open problem. If the two indices $\underline{\mathcal{H}}$ and $\overline{\mathcal{H}}$ are both equal to some $\mathcal{H}\in(0,1)$, methods based on the wavelet decomposition or on discrete filtering (which has several similarities with the wavelet decomposition method) have proved to be often very efficient. The reader is referred to Flandrin (see~\cite{Flan92}), Stoev et al.\ (see~\cite{Sto06}) and the references therein for more informations on the wavelet-based methods and to Kent and Wood (see~\cite{KW97}), Istas and Lang (see~\cite{IL97}) and Coeurjolly (see~\cite{Coeur01,Coeur08}) for more informations about quadratic variations-based methods.

This paper is a first step in the estimation of the two indices $\underline{\mathcal{H}}$ and $\overline{\mathcal{H}}$ in the general case. For this purpose, we investigate the relationship between these two H\"older indices and the wavelet decomposition of a function. The answer is well-known for the index $\underline{\mathcal{H}}$ (see~\cite{Mey90} and theorem~\ref{th:wav-hol} below). The main result of this paper is a characterization of the index $\overline{\mathcal{H}}$, called the upper H\"older exponent, by means of wavelets (see theorem~\ref{th:unifirr} and corollary~\ref{cor:uphold}). Therefore, the results of the present paper should pave the way to the estimation of the indices $\underline{\mathcal{H}}$ and $\overline{\mathcal{H}}$ using wavelet methods.

This paper is organized as follows. In section~\ref{sec:def}, we briefly recall the different concepts for uniform regularity and irregularity. Section~\ref{sec:wav} is devoted to the statement of our main results about the characterization of uniform irregularity by means of wavelets. Finally, section~\ref{sec:proof} contains the proofs of the results stated in section~\ref{sec:wav}.

\section{Upper and lower global H\"older indices}\label{sec:def}
In this section we first give the usual definition of global H\"older index, denoted here lower global H\"older index in order to make a distinction with the upper global H\"older index, which will be introduced afterward.

The definitions rely on the finite differences. For a function $f:\R^d\to\R$ and $x,h\in\R^d$, the first order difference of $f$ is
\[
 \Delta_h^1 f(x)= f(x+h)-f(x).
\]
The difference of order $M$, where $M$ is an integer greater than $2$, is defined by
\[
 \Delta_h^M f(x)= \Delta_h^{M-1} \Delta_h^1 f (x).
\]
Given $\alpha>0$, $\floor{\alpha}$ will denote the greatest integer lower than $\alpha$,
\[
 \floor{\alpha}=\max\{j\in\N\cup\{0\} : j\le \alpha\}.
\]
Throughout this paper, $M$ will designate the integer $M=\floor{\alpha}+1$ and we associate to a bounded function $f:\R^d\to \R$ its $M$-modulus of smoothness $\omega_f^M$:
\[
 \omega_f^M:r\mapsto \sup_{|h|\le r}\sup_{x\in \R^d} |\Delta^M_h f(x)|
\]

\subsection{The lower global H\"older index}\label{sec:lwhi}
Let us recall the well-known notion of lower global H\"older index, usually called global H\"older index or uniform H\"older index.
\BDef
Let $\alpha>0$ and $\beta\in\R$. The bounded function $f$ belongs to $C^{\alpha}_{\beta}(\R^d)$, if there exist $C,R>0$ such that
\begin{equation}\label{eq:def Cab}
 \omega_f^M(r) \le C r^{\alpha}|\log r|^\beta ,
\end{equation}
for any $r\le R$. If $\beta=0$, the space $C^{\alpha}_{0}(\R^d)$ is simply denoted $C^{\alpha}(\R^d)$.

A function $f$ is said to be uniformly H\"{o}lderian if for some $\alpha>0$, $f\in C^{\alpha}(\R^d)$.
\EDef
The above definition leads to a notion of global regularity.
\BDef
The lower global H\"{o}lder exponent of a uniformly H\"{o}lderian function $f$ is defined as
\[
\underline{\mathcal{H}}_f=\sup\{\alpha>0,f\in C^{\alpha}(\R^d)\}.
\]
\EDef

\subsection{The upper global H\"older index}
The irregularity of a function can be studied through the notion of upper global H\"older index. The idea is to reverse inequality~(\ref{eq:def Cab}).

\BDef
Let $f:\R^d\to\R$ be a bounded function, $\alpha\ge 0$ and $\beta\in\R$; $f\in UI^\alpha_{\beta}(\R^d)$ if there exist $C,R>0$ such that
\begin{equation}\label{eq:unifirr}
 \omega_f^M(r) \ge C r^\alpha|\log r|^{\beta}
\end{equation}
for any $r\le R$. If $\beta=0$, the set $UI^\alpha_{0}(\R^d)$ is simply denoted $UI^\alpha(\R^d)$. A function belonging to
$UI^\alpha(\R^d)$ is said to be uniformly irregular with exponent $\alpha$.
\EDef
\BDef
The upper global H\"older exponent (or uniform irregularity exponent) of a bounded function $f$ is
\[
 \overline{\mathcal{H}}_f=\inf\{\alpha:f\in UI^\alpha(\R^d)\}.
\]
\EDef

Let us remark that the statement (\ref{eq:unifirr}) is not a negation of the property $f\in C^\alpha(\R^d)$. Indeed $f$ does not belong to $C^\alpha(\R^d)$ if for any $C>0$, there exists a decreasing sequence $(r_n)_n$ (depending on $C$) converging to $0$ for which
\[
 \omega_f^M(r_n) \ge C r_n^\alpha.
\]
We are thus naturally led to the following definition.
\BDef Let $f:\R^d\to\R$ be a bounded function, $\alpha\ge 0$, $\beta\in\R$; $f\in C^\alpha_{w,\beta}(\R^d)$ if $f\notin UI^\alpha_{\beta}(\R^d)$, i.e.~for any $C>0$ there exists a decreasing sequence $(r_n)_n$ converging to $0$ such that
\[
 \omega_f^M(r_n) \le C r_n^\alpha|\log r_n|^\beta,
\]
for any $n\in\N$. In the case where $\beta=0$, the set $C^\alpha_{w,0}(\R^d)$ is denoted $C^\alpha_{w}(\R^d)$. A function belonging to $C^\alpha_{w}(\R^d)$ is said to be weakly uniformly H\"olderian with exponent $\alpha$.
\EDef
Roughly speaking, a function is weakly uniformly H\"olderian with exponent $\alpha$ if for any $C>0$, one can bound the $M$-modulus of smoothness $\omega^M_f$ of $f$ over $\R^d$ by $\theta(r_n)=C r_n^\alpha |\log r_n|^\beta$ for a remarkable decreasing sequence $(r_n)_n$ of scales, whereas for an H\"{o}lderian function, the $M$-modulus of smoothness of $f$ over $\R^d$ has to be bounded at each scale $r>0$ by $\theta(r)$, for some $C>0$.

\section{A wavelet criterium for uniform irregularity}\label{sec:wav}
In this section we claim that both the lower and upper index of a bounded function can be characterized by means of wavelets.
\subsection{The discrete wavelet transform}
Let us briefly recall some definitions and notations (for more precisions, see e.g.\ \cite{Dau88,Mey90,Mal98}). Under some general assumptions, there exists a function $\phi$ and $2^d-1$
functions $(\psi^{(i)})_{1\le i<2^d}$, called wavelets, such that $\{\phi(x-k)\}_{k\in\Z^d}\cup\{\psi^{(i)}(2^j x-k):1\le i<2^d, k\in \Z^d, j\in\Z \}$ form an orthogonal basis of
$L^2(\R^d)$. Any function $f\in L^2(\R^d)$ can be decomposed as follows,
\[
f(x)=\sum_{k\in \Z^d} C_k \phi(x-k) + \sum_{j=1}^{+\infty}
\sum_{k\in\Z^d} \sum_{1\le i<2^d} c^{(i)}_{j,k} \psi^{(i)}(2^j
x-k),
\]
where
\[
c^{(i)}_{j,k}=2^{dj}\int_{\R^d}f(x) \psi^{(i)}(2^jx-k)\, dx,
\]
and
\[
C_k=\int_{\R^d} f(x) \phi(x-k)\, dx.
\]
Let us remark that we do not choose the $L^2(\R^d)$ normalization for the wavelets, but rather an $L^\infty$ normalization, which is better fitted to the study of the H\"olderian regularity. Hereafter, the wavelets are always supposed to belong to $C^\gamma(\R^d)$ with $\gamma$ sufficiently large (we require at least $\gamma>\alpha$) and the functions $\{\partial^{s}\phi\}_{|s|\le \gamma}$, $\{\partial^{s}\psi^{(i)}\}_{|s|\le \gamma}$ are assumed to have fast decay. Furthermore, in $\R^d$ we will use the tensor product wavelet basis (see \cite{Mey90,Dau92} and section~\ref{sec:proofpro}).

A dyadic cube of scale $j$ is a cube of the form
\[
\lambda=\left[\frac{k_1}{2^j},\frac{k_1+1}{2^j}\right)\times
\cdots \times \left[\frac{k_d}{2^j},\frac{k_d+1}{2^j}\right),
\]
where $k=(k_1,\ldots,k_d)\in \Z^d$. From now on, wavelets and wavelet coefficients will be indexed with dyadic cubes $\lambda$. Since $i$ takes $2^d-1$ values, we can assume that it takes values in $\{0,1\}^d\setminus\{(0,\ldots,0)\}$; we will use the following notations:
\begin{itemize}
\item
$\lambda=\lambda(i,j,k)=\frac{k}{2^j}+\frac{i}{2^{j+1}}+[0,\frac{1}{2^{j+1}})^d$,
\item $c_\lambda=c^{(i)}_{j,k}$, \item
$\psi_\lambda=\psi^{(i)}_{j,k}=\psi^{(i)}(2^j\cdot -k)$.
\end{itemize}

To state our wavelet criteria, we will use the following notation: for any $j\geq 0$, we set
\[
\|c_{j,\cdot}^{(\cdot)}\|_{\infty}=\sup_{i\in\{0,1\}^d\setminus\{(0,\ldots,0)\}}\sup_{k\in\Z^d}|c_{j,k}^{(i)}|\;.
\]

\subsection{Wavelets and usual uniform regularity}
The characterization of the lower global H\"older index in terms of wavelet coefficients is well-known.

The uniform H\"olderian regularity of a function is closely related to the decay rate of its wavelet coefficients. Let us recall the following result (see \cite{Mey90}).
\BThe\label{th:wav-hol} Let $\alpha>0$ such that $\alpha\not\in\N$. We have $f\in C^{\alpha}(\R^d)$ if and only if there exists $C>0$ such that
\begin{equation}\label{eq:h:cond}
\left\{
\begin{array}{l}
\forall k\in\Z^d,\,|C_k|\leq
C\;,\\
\forall j\geq 0,\;\|c_{j,\cdot}^{(\cdot)}\|_{\infty}\leq
C2^{-j\alpha}\;.
\end{array}\right.
\end{equation}
\EThe

This theorem yields a wavelet characterization of the lower H\"older index of a uniformly H\"olderian function.
\BCor Assume that $f$ is a uniformly H\"olderian function; we have
\[
\underline{\mathcal{H}}_f=\liminf_{j\to \infty}\frac{\log_2 \|c_{j,\cdot}^{(\cdot)}\|_\infty}{-j} .
\]
\ECor

\subsection{Wavelets and uniform irregularity}
In this section, we aim at characterizing the uniform irregularity of a bounded function in terms of wavelets.

The main result is the following theorem.
\BThe\label{th:unifirr}
Let $\alpha>0$ and $f$ be a bounded function on $\R^d$. If there exists $C>0$ such that for any integer $j\geq 0$,
\begin{equation}\label{eq:mino}
\max\left(\;\sup_{\ell\ge j} \|c_{\ell,\cdot}^{(\cdot)}\|_{\infty}
 ,\;2^{-jM}\sup_{\ell\le  j} (2^{\ell M}\|c_{\ell,\cdot}^{(\cdot)}\|_{\infty})\right)\ge C 2^{-j\alpha},
\end{equation}
then $f\in UI^{\alpha}(\R^{d})$.

Conversely, if $f$ is uniformly H\"{o}lder and if for any $\beta>1$, $f$ belongs to $UI^{\alpha}_{\beta}(\R^{d})$, then there exists $C>0$ such that relation~(\ref{eq:mino}) holds for any $j\geq 0$.
\EThe
Let us make some remarks.
\BRem
Unlike the case of usual uniform H\"olderian regularity, the case where $\alpha$ is a natural number is not a specific one.
\ERem
\BRem
The assumptions of Theorem~\ref{th:unifirr} are indeed optimal. See Section~\ref{sec:cex} in Appendix for more details.
\ERem
\BRem
The condition
\[
\|c_{j,\cdot}^{(\cdot)}\|_{\infty}\ge C 2^{-j\alpha},
\]
for some $C>0$ and any $j\ge 0$ is a sufficient (but not necessary) condition for uniform irregularity. In the general case,
\[
\overline{\mathcal{H}}_f\neq \limsup_{j\to
+\infty}\frac{\log_2 \|c_{j,\cdot}^{(\cdot)}\|_{\infty} }{-j}
\;.
\]
\ERem

Following theorem~\ref{th:wav-hol}, a bounded function $f$ is not uniformly H\"olderian with exponent $\alpha$, i.e.~its $M$-modulus of smoothness is bounded from below by $\theta(r_n)$ for some specific decreasing sequence $(r_n)$ converging to
$0$, if and only if a similar property holds for its wavelet coefficients. The situation is completely different concerning uniform irregularity: the value of the $M$-modulus of smoothness at $r=2^{-j}$ is influenced by the wavelet coefficients at scales below and above the scale $2^{-j}$. The $M$-modulus of smoothness of $f$ can be large at $r=2^{-j}$ for any $j\in\N$ (even if for some scales $j$, the coefficients $(c_{j,k}^{(i)})$ are small or even vanish) provided that for any $j\in\N$, at a controlled distance of the scale $2^{-j}$, there exists some large wavelet coefficients. Such a behavior is met with the lacunary fractional Brownian motion, which admits some vanishing wavelet coefficients but that is almost surely locally uniformly irregular (see~\cite{Clau10} for more details).

Theorem~\ref{th:unifirr} leads to a wavelet characterization of the upper H\"older exponent.
\BCor\label{cor:uphold}
If $f$ is a uniformly H\"olderian function, then
\[
\overline{\mathcal{H}}_f=\limsup_{j\to\infty}
\frac{\log_2 \max\left(\sup_{\ell\ge j}\|c_{\ell,\cdot}^{(\cdot)}\|_{\infty},
 2^{-j M}\sup_{\ell\le j} (2^{\ell M}\|c_{\ell,\cdot}^{(\cdot)}\|_{\infty})\right)}{-j} .
\]
\ECor

\section{Proof of Theorem~\ref{th:unifirr}}\label{sec:proof}
We will show that theorem~\ref{th:unifirr} comes from the following wavelet characterization (up to a logarithmic term) of the weak uniform H\"olderian regularity.
\BPro\label{pro:wavcarWR}
Let $\alpha>0$;
\begin{enumerate}
\item if $f\in C_w^\alpha(\R^d)$ then, for any $C>0$, there exists a strictly increasing sequence of integers $(j_n)_{n\in \N}$ such that for any $n\ge 0$ and any $j\in \{j_n,\ldots,j_{n+1}-1\}$,
\begin{equation}\label{eq:wav-carWR}
\sup_{|\lambda|=2^{-j}} |c_{\lambda}| \le C'C\inf(2^{-j_n \alpha},2^{(M-\alpha)j_{n+1}} 2^{-jM}),
\end{equation}
for some $C'>0$ depending only on the chosen wavelet basis.
\item Conversely, if $f$ is uniformly H\"olderian and if for any $C>0$, there exists a strictly increasing sequence of integers $(j_n)_{n\in \N}$ such that~(\ref{eq:wav-carWR}) holds then $f\in C^\alpha_{w,\beta}(\R^d)$ for any $\beta>1$.
\end{enumerate}
\EPro

\subsection{A reformulation of the property $f\in C_w^\alpha(\R^d)$}\label{sec:fineproperty}
To prove Proposition~\ref{pro:wavcarWR}, we first need to reformulate in a more appropriate way the property $f\in C_w^\alpha(\R^d)$.

Since modulus of smoothness $\omega_f^M$ is a non-decreasing function, $f\in C_w^\alpha(\R^d)$ if and only if, for any $C>0$, there exists an increasing sequence of integers $(j_n)_{n\in\N}$ such that for any $r\in (2^{-j_{n+1}},2^{-j_n}]$ ($n\in\N$),
\begin{equation}\label{eq:ubcoarse}
\omega_f^M(r)=\sup_{|h|\le r}\sup_{x\in \R^d} |\Delta^M_h f(x)| \le C 2^{-j_n \alpha}.
\end{equation}
Hence, $f$ belongs to $C_w^\alpha(\R^d)$ if and only if the piecewise constant function $\Theta$ defined as
\[
\Theta = C \sum_{n\in\N} 2^{-j_n \alpha} \chi_{(2^{-j_{n+1}},2^{-j_n}]},
\]
where $\chi_A$ denotes the characteristic function of the set $A$ is an upper bound of the $M$-modulus of smoothness $\omega^M_f$ of $f$.

This characterization of the weak uniform regularity is not convenient to deal with, since
\[
\limsup_{r\to 0} \frac{\Theta(2r)}{\Theta(r)}
\]
may be infinite. To overcome this problem, in the next proposition we will reformulate the property $f\in C_w^\alpha(\R^d)$, giving a finer upper bound of $\omega_f^M$. To this end, let us remark that there is a link between the finite differences of $f$ at different scales.
\BPro\label{pro:reformulationWR}
The bounded function $f$ belongs to $C_w^\alpha(\R^d)$ if and only if for any $C>0$, there exists a strictly increasing sequence of integers $(j_n)_{n\in\N}$ such that for any $j\in \{j_n,\ldots,j_{n+1}-1\}$,
\begin{equation}\label{eq:ubfine}
\sup_{|h|\le 2^{-j}}\sup_{x\in \R^d}| \Delta^M_h f(x)| \le C \inf(2^{-j_n \alpha},2^{M(j_{n+1}-j)} 2^{-j_{n+1}\alpha}).
\end{equation}
\EPro
\BProof
Let us first assume that (\ref{eq:ubcoarse}) holds. The following relation (given in \cite{PP87} for example),
\[
\omega_f^M (2r)=\sup_{|h|\le 2r} \sup_{x\in \R^d} |\Delta^M_h f(x)| \le  2^M \sup_{|h|\le r}\sup_{x\in \R^d} |\Delta^M_h f(x)|= 2^M \omega_f^M (r)
\]
and equation~(\ref{eq:ubcoarse}) imply that for any $j\in \{j_n,\ldots,j_{n+1}-1\}$,
\begin{eqnarray*}
\omega^M_f (2^{-j}) &=& \omega^M_f (2^{j_{n+1}-j}2^{-j_{n+1}}) \\
 &\le& 2^{M(j_{n+1}-j)}\omega^M_f (2^{-j_{n+1}})\le C 2^{M(j_{n+1}-j)} 2^{-j_{n+1}\alpha}.
\end{eqnarray*}
Hence, relation~(\ref{eq:ubfine}) holds. The converse assertion is obvious.
\EProof

Let us now remark that the piecewise function $\theta$ defined (on $(0,2^{-j_1}]$) as
\begin{equation}\label{eq:theta}
\theta (r)= \sum_{n\in\N} \inf(2^{-j_n \alpha}, 2^{j_{n+1} (M-\alpha)}r^M) \chi_{(2^{-j_{n+1}}, 2^{-j_n}]}(r)
\end{equation}
is a continuous function. Furthermore it satisfies additional interesting properties summed up in the following proposition.
\BPro\label{pro:modcont}
Let $\alpha>0$ and $(j_n)_{n\in\N}$ be an increasing sequence of integers. Let $\theta$ be defined by equality~(\ref{eq:theta}). The function $\theta$ obeys the following properties:
\begin{enumerate}
\item $\theta$ is a modulus of continuity, that is a non decreasing continuous function satisfying
\begin{equation}\label{eq:modcont}
\limsup_{r\to  0} \frac{\theta(2r)}{\theta(r)} <\infty,
\end{equation}
 \item for any $\beta>1$ and for any $J$ sufficiently large, the following relations are satisfied:
\begin{equation}\label{eq:modcontFaibleUn}
 \sum_{j=j_1}^J 2^{Mj}\theta(2^{-j})\le C J2^{MJ}\theta(2^{-J})
\end{equation}
\begin{equation}\label{eq:modcontFaibleDeux}
 \sum_{j\ge J} \frac{\theta(2^{-j})|\log \theta(2^{-j})|^\beta}{j^\beta} \le C J^{\beta}\theta(2^{-J}),
\end{equation}
\begin{equation}\label{eq:modcontFaibleTrois}
 2^{-Mj}=o(\theta(2^{-j}))\mbox{ as }j\to \infty.
\end{equation}
\end{enumerate}
\EPro
\BProof
We first prove that $\theta$ is a modulus of continuity by showing that
\begin{equation}\label{eq:modcontbis}
\theta(2r)\le 2^M \theta(r).
\end{equation}
Assume that there exists some $n\in\N$ such that
\[
 2^{-j_{n+1}}\le r\le 2^{-j_n-1}.
\]
Since $2^{-j_{n+1}+1}\le 2r\le 2^{-j_n}$, one has
\[
\theta(2r)=\inf(2^{-j_n \alpha},2^{j_{n+1}(M-\alpha)}(2r)^M)\leq 2^M \theta(r).
\]
On the other hand, if for some $n\in\N$, one has
\[
 2^{-j_n-1}\le r\le 2^{-j_n},
\]
then $2^{-j_n}\le 2r \le 2^{-j_n+1}$ and thus
\begin{eqnarray*}
 \theta(2r) &=& \inf(2^{-j_{n-1}\alpha}, 2^{j_n (M-\alpha)}(2r)^M)\\
 &\le& 2^M (2^{j_n}r)^M 2^{-j_n \alpha} = 2^M 2^{j_n (M-\alpha)}r^M.
\end{eqnarray*}
Since $M-\alpha>0$, one has
\[
 2^M 2^{j_n (M-\alpha)}r^M \le 2^M 2^{j_{n+1} (M-\alpha)}r^M.
\]
Moreover, since $r\le 2^{-j_n}$,
\[
 2^M (2^{j_n}r)^M 2^{-j_n \alpha}\le 2^M 2^{-j_n\alpha},
\]
hence,
\[
 \theta(2r)\leq
 2^M \inf(2^{-j_n \alpha},2^{j_{n+1}(M-\alpha)}r^{M})\;.
\]
In any case, relation~(\ref{eq:modcontbis}) holds, which directly implies~(\ref{eq:modcont}).

Let us now prove the second part of Proposition~\ref{pro:modcont}. Let $J\in\N$ and $n_0\in\N$ such that $j_{n_0}\le J\le j_{n_0+1}-1$. Let us first show that property~(\ref{eq:modcontFaibleUn}) is satisfied. By definition, we have
\begin{eqnarray*}
 \sum_{j=j_1}^J 2^{Mj}\theta(2^{-j}) &=&
 \sum_{n=0}^{n_0-1}\sum_{j=j_n}^{j_{n+1}-1} 2^{Mj} \inf(2^{-j_n \alpha},2^{j_{n+1} (M-\alpha)} 2^{-jM}) \\
 && +\sum_{j=j_{n_0}}^{J-1} 2^{Mj}\inf(2^{-j_{n_0}\alpha},2^{j_{n_0+1} (M-\alpha)} 2^{-jM}).
\end{eqnarray*}
Therefore,
\[
 \sum_{j=j_1}^J 2^{Mj} \theta(2^{-j})\le \sum_{n=0}^{n_0-1} j_{n+1} 2^{j_{n+1}(M-\alpha)}
 + J \inf(2^{MJ} 2^{-j_{n_0}\alpha},2^{j_{n_0+1} (M-\alpha)}),
\]
that is
\begin{eqnarray*}
 \sum_{j=j_1}^J 2^{Mj} \theta(2^{-j}) &\le&
 j_{n_0} 2^{j_{n_0} (M-\alpha)} +J\inf(2^{MJ} 2^{-j_{n_0}\alpha},2^{j_{n_0+1} (M-\alpha)}) \\
 &\leq& 2J \inf(2^{MJ} 2^{-j_{n_0} \alpha}, 2^{j_{n_0+1} (M-\alpha)}),
\end{eqnarray*}
which shows that property~(\ref{eq:modcontFaibleUn}) holds.

We now check inequality~(\ref{eq:modcontFaibleDeux}). Since
\[
 \theta(2^{-j})\leq 2^{-j_n\alpha}
\]
for any $n\ge n_0$ and any $j\in\{j_n,\ldots,j_{n+1}-1\}$,
we have
\begin{eqnarray*}
 \sum_{j=J}^\infty \frac{\theta(2^{-j}) |\log \theta(2^{-j})|^\beta} {j^{\beta}}
 &\le& \sum_{j=J}^{j_{n_0+1}-1} \frac{\theta(2^{-j}) |\log\theta(2^{-j})|^\beta}{j^\beta} \\
 && +\sum_{n=n_0+1}^{\infty} \sum_{j=j_n}^{j_{n+1}-1} \frac{2^{-j_n \alpha} |\log(2^{-j_n\alpha})|^\beta}{j^\beta}\\
 &=& \sum_{j=J}^{j_{n_0+1}-1} \frac{\theta(2^{-j})|\log\theta(2^{-j})|^\beta}{j^\beta} \\
 && + C \sum_{n=n_0+1}^{\infty} 2^{-j_n \alpha}j_n^\beta \sum_{j=j_n}^{j_{n+1}-1} \frac{1}{j^\beta}.
\end{eqnarray*}
Using equality~(\ref{eq:theta}), we get
\begin{eqnarray}\label{eq:modcontproof1}
 \sum_{j=J}^{\infty} \frac{\theta(2^{-j}) |\log\theta(2^{-j})|^\beta}{j^\beta} &\le&
 C \sum_{j=J}^{j_{n_0+1}-1} \frac{\inf(j_{n_0}^\beta 2^{-j_{n_0}\alpha} ,j^\beta 2^{j_{n_0+1} (M-\alpha)}2^{-jM})}{j^\beta} \nonumber \\
 &+& C \sum_{n=n_0 +1}^{\infty} j_n^\beta 2^{-j_n \alpha} \sum_{j=j_n}^{j_{n+1}-1}\frac{1}{j^\beta}.
\end{eqnarray}
Moreover, since
\[
 \sum_{n=n_0+1}^{\infty} j_n^\beta 2^{-j_n \alpha}\sum_{j=j_n}^{j_{n+1}-1} \frac{1}{j^\beta} \le
 \sum_{n=n_0+1}^{\infty} j_n 2^{-j_n \alpha}\le
 j_{n_0+1} 2^{-j_{n_0+1}\alpha},
\]
inequality~(\ref{eq:modcontproof1})
 yields
\begin{eqnarray*}
 \sum_{j=J}^{\infty} \frac{\theta(2^{-j}) |\log\theta(2^{-j})|^\beta}{j^\beta} &\le&
 C \sum_{j=J}^{j_{n_0+1}-1} \frac{\inf(j_{n_0}^\beta 2^{-j_{n_0} \alpha}, j^\beta 2^{j_{n_0+1} (M-\alpha)} 2^{-jM})}{j^\beta}\\
 &&+ C j_{n_0+1} 2^{-j_{n_0+1}\alpha} \\
 &\le& C' (\inf(j_{n_0}2^{-j_{n_0 \alpha}}, 2^{j_{n_0+1}(M-\alpha)}2^{-JM}) \\ && +j_{n_0+1}2^{-j_{n_0+1}\alpha})\\
&\le& C' J^\beta \theta(2^{-J}).
\end{eqnarray*}

Since $M>\alpha$, relation~(\ref{eq:modcontFaibleTrois}) is straightforward.
\EProof
\BRem
The concept of modulus of continuity has been used in~\cite{JafMey96} to deal with a more general notion of uniform H\"olderian regularity than the usual one, induced by the H\"older spaces. For a given $M$ and a given modulus of continuity $\theta$, a wavelet characterization of the property
\begin{equation}\label{eq:modcontReg}
\omega_f^M (r)\le C \theta(r)
\end{equation}
for any $r\ge 0$ is provided under the two following assumptions on $\theta$:
for any $J\ge 0$,
\begin{equation}\label{eq:modcontFortUn}
 \sum_{j=0}^J 2^{jM} \theta(2^{-j})\le C' 2^{JM}\theta(2^{-J})
\end{equation}
and
\begin{equation}\label{eq:modcontFortDeux}
 \sum_{j=J}^{\infty} 2^{j(M-1)} \theta(2^{-j})\le C' 2^{J(M-1)} \theta(2^{-J}).
\end{equation}
Properties~(\ref{eq:modcontFortUn}) and (\ref{eq:modcontFortDeux}) are much stronger than properties~(\ref{eq:modcontFaibleUn}), (\ref{eq:modcontFaibleDeux}) and (\ref{eq:modcontFaibleTrois}), which concern the weak uniform regularity of a function $f$.
\ERem

\subsection{Proof of Proposition~\ref{pro:wavcarWR}}\label{sec:proofpro}
We shall split the proof into two parts.

\BPro
Let $\alpha>0$; if $f\in C_w^\alpha(\R^d)$ then, for any $C>0$, there exists a strictly increasing sequence of integers $(j_n)_{n\in \N}$ such that for any $n\ge 0$ and any $j\in \{j_n,\ldots,j_{n+1}-1\}$,
\[
 \sup_{|\lambda|=2^{-j}} |c_{\lambda}| \le C'C \theta(2^{-j}),
\]
for some $C'>0$ depending only on the chosen wavelet basis, where $\theta$ is the function defined by equality~(\ref{eq:theta}).
\EPro
\BProof
Assume that $f$ belongs to $C^\alpha_w(\R^d)$ and let $C>0$. By proposition~\ref{pro:reformulationWR}, we have for any $r$ sufficiently small,
\begin{equation}\label{eq:proofprop1a}
\omega^M_f (r)\leq C \theta(r).
\end{equation}

If $d=1$, let us recall (see \cite{Jaf04}) that if the wavelet basis belongs to $C^M(\R^d)$ then there exists a function $\Psi_{M}$ with fast decay and such that $\psi=\Delta_{\frac{1}{2}}^{M}\Psi_{M}$.
In dimension $d>1$, we use the tensor product wavelet basis:
\[
\psi^{(i)}(x)=\Psi^{(1)}(x_{1})\cdots\Psi^{(d)}(x_{d}),
\]
where for all $i$, $\Psi^{(i)}$ are either $\psi$ or $\phi$ but at least one of them must equal $\psi$. For example, assume that $\Psi^{(1)}=\psi$. Then, for any $i\in \{1,\ldots,2^{d}-1\}$, any $j\ge 0$ and any $k\in \Z^d$,
\[
 c_{j,k}^{(i)} = 2^{jd} \int_{\R^d} f(x) \Psi^{(1)}(2^j x_1 -k_1) \cdots \Psi^{(d)}(2^j x_d -k_d) \,dx.
\]
We thus have
\begin{eqnarray*}
 c_{j,k}^{(i)} &=& 2^{jd} \int_{\R^d} f(x) \Delta_{1/2}^M \Psi_M (2^j x_1 -k_1) \cdots \Psi^{(d)} (2^j x_d -k_d) \,dx \\
 &=& 2^{jd} \int_{\R^d} \Delta_{1/2^{j+1}e_1}^M f (x) \Psi_M (2^j x_1 -k_1) \cdots \Psi^{(d)} (2^j x_d -k_d) \, dx,
\end{eqnarray*}
with $e_1=(1,0,\cdots,0)$ and therefore
\[
 |c_{j,k}^{(i)}|\le 2^{jd}\int_{\R^d} |\Delta_{1/2^{j+1}e_1}^Mf (x)| |\Psi_M(2^j x_1 -k_1)\cdots \Psi^{(d)} (2^j x_d -k_d)| \,dx.
\]
We thus get, using inequality~(\ref{eq:proofprop1a}),
\[
 |c_{j,k}^{(i)}|\le C 2^{jd} \theta(2^{-(j+1)}) \int_{\R^d} |\Psi_M (2^j x_1 -k_1)\cdots \Psi^{(d)} (2^j x_d -k_d)| \,dx.
\]
Setting $y=2^j x-k$ in the last integral, we obtain
\[
 2^{jd} \int_{\R^d} |\Psi_M (2^j x_1 -k_1) \cdots \Psi^{(d)} (2^j x_d -k_d)| \,dx
 = \|\Psi_M \otimes \cdots \otimes \Psi^{(d)} \|_{L^1(\R^d)}.
\]
Since $\theta$ is a non-decreasing function, we can write
\[
 |c_{j,k}^{(i)}| \le C \theta(2^{-j}) \|\Psi_{M}\|_{L^1(\R^d)},
\]
which ends the proof.
\EProof

From now on in this section, we suppose that $f$ is uniformly H\"olderian and that property~(\ref{eq:wav-carWR}) is satisfied. For the second part of the proof, we need to introduce the following notations:
\begin{equation}\label{Eqfj}
 f_{-1} (x) = \sum_{k\in \Z^d} C_k \varphi(x-k) ,\quad
 f_{j}(x)= \sum_{i=1}^{2^d-1} \sum_{k\in\Z^d} c_{j,k}^{(i)} \psi(2^jx-k),
\end{equation}
with $j\ge 0$. Since $f$ is uniformly H\"olderian, $f_j$, as defined by equality~(\ref{Eqfj}), converges uniformly on any compact to a limit which has the same regularity as the wavelets. Furthermore $\sum_{j\ge -1}f_j (x)$ converges uniformly on any compact. The proof is based on the following lemma which provides an upper bound for $\|\partial^\gamma f_j (x)\|_{L^\infty(\R^d)}$, for any $|\gamma|\le M$.
\BLem\label{lem:ubfjm}
Let $m\in\{0,\ldots,M\}$; there exists some $C'>0$ depending only on $m$ and on the chosen wavelet basis such that for any $\gamma\in\N^d$ satisfying $|\gamma|=m$ and for $j$ sufficiently large,
\[
\|\partial^\gamma f_j (x)\|_{L^\infty(\R^d)} \le C' C 2^{jm} \inf(\theta(2^{-j}), \frac{\theta(2^{-j})|\log\theta(2^{-j})|^\beta}{j^\beta}),
\]
where $\theta$ is the function defined by equality~(\ref{eq:theta}).
\ELem
\BProof
Since $f$ satisfies Property~(\ref{eq:wav-carWR}), one has
\begin{equation}\label{eq:proofprop1c}
 |c_{j,k}^{(i)}|\le C \theta(2^{-j}),
\end{equation}
for $j$ sufficiently large. Furthermore, since $f$ is uniformly H\"olderian,
\begin{equation}\label{eq:proofprop1d}
 \big|\log |c_{j,k}^{(i)}|\big|\ge C' j,
\end{equation}
for some $C'>0$ and $j$ sufficiently large. Now, using the trivial relation
\[
 |c_{j,k}^{(i)}| =\inf\left(|c_{j,k}^{(i)}|, \frac{|c_{j,k}^{(i)}|\big|\log |c_{j,k}^{(i)}|\big|^\beta}{\big|\log|c_{j,k}^{(i)}|\big|^\beta}\right),
\]
inequalities~(\ref{eq:proofprop1c}) and~(\ref{eq:proofprop1d}) leads to
\[
 |c_{j,k}^{(i)}|\le \inf\left(\theta(2^{-j}), \frac{\theta(2^{-j}) |\log \theta(2^{-j})|^\beta}{j^\beta}\right).
\]
Therefore, for any integer $p>d$,
\begin{eqnarray*}
|\partial^\alpha f_j(x)| &=& |\sum_{i=1}^{2^d -1}\sum_{k\in\Z^d} c_{j,k}^{(i)} 2^{jm} \partial^\alpha \psi^{(i)} (2^{j}x-k)| \\
&\le& C'C 2^{jm} \sum_{i=1}^{2^d -1} \sum_{k\in\Z^d} \frac{\inf (\theta(2^{-j}),\frac{\theta(2^{-j}) |\log\theta(2^{-j})|^\beta}{j^\beta})}{(1+|2^j x-k|)^p},
\end{eqnarray*}
using the fast decay of the wavelets. The use of the classical bound
\[
\sup_{x\in\R^d} \sum_{k\in\Z^d} \frac{1}{(1+|2^j x-k|)^p}<\infty
\]
ends the proof of this lemma.
\EProof
\BPro
Let $\alpha>0$; if $f$ is uniformly H\"olderian and if for any $C>0$, there exists a strictly increasing sequence of integers $(j_n)_{n\in \N}$ such that~(\ref{eq:wav-carWR}) holds, let $h\in \R^d$ and define $J=\sup\{j_n: |h|< 2^{-j_n}\}$. We have, for $h$ sufficiently small,
\begin{equation}\label{eq:proofprop1b}
|\Delta^M_h f(x)|\le C' J^\beta\theta(2^{-J}),
\end{equation}
where $\theta$ is the function defined by equality~(\ref{eq:theta}).
\EPro
\BProof
Let us set
\[
g_1=\sum_{j=-1}^{j_1 -1} f_j(x),\quad g_2=\sum_{j=j_1}^{J-1} f_j, \quad\mbox{and}\quad g_3=\sum_{j=J}^\infty \Delta^M_h f_j(x).
\]
For any $j\ge -1$, $f_j$ has the same regularity as the wavelets and so does $g_1$. Therefore, we can suppose that $g_1$ belongs to $C^\eta (\R^d)$ with $M<\eta\notin \N$ and for any $r>0$,
\[
\omega^M_{g_1} (r)\le C' r^M,
\]
(see e.g.~\cite{kra83}). Using relation~(\ref{eq:modcontFaibleTrois}), we get that inequality~(\ref{eq:proofprop1b}) holds for $f=g_1$.

Let us now consider the case $f=g_2$. Lemma~\ref{lem:ubfjm} with $m=M$ leads to the inequality
\[
|\partial^{\gamma} f_j (x)|\le C' C 2^{jM} \theta(2^{-j})
\]
for any $\gamma$ such that $|\gamma|=M$ and for any $j_1\le j\leq J-1$. Furthermore, for any $j$, $f_j\in C^\eta(\R^d)$ which can be considered as a subset of the homogeneous H\"older space $\dot{C}^\eta(\R^d)$ (see e.g.~\cite{Mey90}). Therefore,
\begin{eqnarray*}
 |\Delta^M_h f_j(x)|
 \le|h|^M \sum_{|\gamma|=M} \|\partial^\gamma f_j\|_{L^\infty(\R^d)},
\end{eqnarray*}
for any $j\geq j_1$. We thus have
\[
 |\sum_{j=j_0}^{J-1} \Delta^M_h f_j(x)|\le
 C'C |h|^M \sum_{j=j_0}^{J-1} 2^{jM} \theta(2^{-j}).
\]
Using relation~(\ref{eq:modcontFaibleUn}), we get
\[
 |\sum_{j=j_0}^{J-1} \Delta^M_h f_j (x)|
 \le C'C |h|^M J 2^{JM} \theta(2^{-J})
 \le C'C J \theta(2^{-J}).
\]
We have thus proved that the function $g_2$ satisfies inequality~(\ref{eq:proofprop1b}).

For $g_3$, let us apply lemma~\ref{lem:ubfjm} with $m=0$ to obtain
\[
 |\sum_{j=J}^{\infty} \Delta^M_h f_j(x)| \le C' C \sum_{j=J}^{\infty} \frac{\theta(2^{-j})|\log\theta(2^{-j})|^\beta}{j^\beta}.
\]
By inequality~(\ref{eq:modcontFaibleDeux}), we have
\[
|\sum_{j=J}^{\infty} \Delta^M_h f_j(x)|\le C'C J^\beta \theta(2^{-J}).
\]

The results concerning $g_1,g_2$ and $g_3$ put together show that the function $f$ satisfies inequality~(\ref{eq:proofprop1b}), which ends the proof.
\EProof

\subsection{Proof of Theorem~\ref{th:unifirr}}\label{sec:proofmain}
We now prove that Theorem~\ref{th:unifirr}, characterizing the uniform irregularity in terms of wavelet coefficients, is the contrapositive of proposition~\ref{pro:wavcarWR}.

We just need the following lemma.
\BLem\label{lem:contrap}
The two following assertions are equivalent:
\begin{enumerate}
\item the wavelet coefficients of $f$ do not satisfy property~(\ref{eq:wav-carWR}),
\item there exists $C'>0$ and an integer $j_0$ such that, for any $j\ge j_0$,
\begin{equation}\label{eq:minowc}
 \max\big(\sup_{\ell\ge j}\sup_{|\lambda|=2^{-\ell}} |c_{\lambda}|
 ,2^{-jM}\sup_{\ell\le j} (2^{\ell M}\sup_{|\lambda|=2^{-\ell}} |c_{\lambda}|)\big)
 \ge C'\theta(2^{-j}).
\end{equation}
\end{enumerate}
\ELem
\BProof
Let us show that property~(\ref{eq:wav-carWR}) is equivalent to the negation of property~(\ref{eq:minowc}). Indeed by definition, the wavelet coefficients of $f$ satisfy property~(\ref{eq:wav-carWR}) if and only if for any $C>0$, there exists an increasing sequence of integers $(j_{n})_{n\in\N}$ such that
\[
 \sup_{i,k} |c_{j,k}^{(i)}| \le
 C \inf(2^{-j_n \alpha},2^{j_{n+1}(M-\alpha)}2^{-jM}),
\]
for any $n\in\N$ and any $j\in \{j_{n},\cdots,j_{n+1}-1\}$. This statement can be reformulated as follows: for any $C>0$, there exists an increasing sequence of integers $(j_n)_{n\in\N}$ such that for any $n\in\N$,
\[
 \sup_{\ell\ge j_n} \sup_{i,k} |c_{\ell,k}^{(i)}| \le C 2^{-j_n\alpha}
\]
and
\[
 \sup_{j_0\le \ell\le j_{n+1}} 2^{\ell M}\sup_{i,k} |c_{\ell,k}^{(i)}| \le C 2^{j_{n+1}(M-\alpha)}.
\]
Let us set
\[
 n_0=\inf\{n\in\N: \sup_{0\le \ell\le j_0} (2^{\ell M}\sup_{i,k}|c_{\ell,k}^{(i)}|)\le C 2^{j_{n+1}(M-\alpha)}\}.
\]
Replacing the sequence $j_n$ by $\ell_n=j_{n+n_0+1}$, property~(\ref{eq:wav-carWR}) is equivalent to the existence, for any $C>0$, of a strictly increasing sequence of integers $(j_n)_{n\in\N}$ such that for any  $n\in\N$,
\[
 \sup_{\ell\ge j_n} \sup_{i,k} |c_{\ell,k}^{(i)}| \le C 2^{-j_n\alpha},
\]
and
\[
 \sup_{\ell\le j_n} 2^{\ell M}\sup_{i,k} |c_{\ell,k}^{(i)}| \le C 2^{j_n (M-\alpha)}.
\]

To conclude, observe that the last property is equivalent to the existence, for any $C>0$ and any $j_0\in \N$, of some $j_1>j_0$ such that
\[
 \sup_{\ell\ge j_1} \sup_{i,k} |c_{\ell,k}^{(i)}| \le C 2^{-j_1\alpha}
\]
and
\[
 \sup_{\ell\le j_1} 2^{\ell M} \sup_{i,k} |c_{\ell,k}^{(i)}| \le C 2^{j_1 (M-\alpha)}.
\]
Since this is the negation of relation~(\ref{eq:minowc}), the lemma is proved.
\EProof

Theorem~\ref{th:unifirr} directly follows from Proposition~\ref{pro:wavcarWR} and Lemma~\ref{lem:contrap}.

\appendix
\section{Optimality of the assumptions of Theorem~\ref{th:unifirr}}\label{sec:cex}
We prove here the optimality of the assumptions of proposition~\ref{pro:wavcarWR} and thus of theorem~\ref{th:unifirr}. To this end we use two counter-examples already introduced in ~\cite{Jaf89}.

\subsection{A uniform irregular function satisfying Property~(\ref{eq:wav-carWR})}
Let $\alpha\in (0,1)$, $\ell_0\in\N$ and define the two following sequences of integers $(j_n)_{n\in\N}$ and $(j_{n,\alpha})_{n\in\N}$ as
\[
\left\{
 \begin{array}{ll}
 j_1=\ell_0,\\
 j_{n+1}=[\frac{1}{1-\alpha}2^{j_n \alpha}-j_n \alpha],& \forall n\geq 1,\\
 j_{n,\alpha}=\floor{2^{j_n \alpha}},& \forall n\geq 1.
 \end{array}
\right.
\]
We aim at proving the following result.
\BPro\label{pro:cex1}
Let us assume that the multiresolution analysis is compactly supported. Let $\varepsilon\in (0,1)$ and $\ell_0$ be such that $\mathrm{supp}(\psi)\subset [-2^{\ell_0},2^{\ell_0}]$. Furthermore, let us assume that $\psi(0)\neq 0$. The function $f$ defined as
\begin{eqnarray*}
 f(x) &=& \sum_{n=0}^{\infty} 2^{-j_n\alpha} \sum_{j=j_n}^{j_{n,\alpha}} \sum_{\ell=j+2}^{j_{n,\alpha}} \ell^{-\varepsilon} \psi\big(2^{\ell}(x-2^{-(j-\ell_0)}) \big)\\
 & & +\sum_{n=0}^{\infty} 2^{j_{n+1}(1-\alpha)}\sum_{j=j_{n,\alpha}+1}^{j_{n+1}-1} \sum_{\ell=j+2}^{j_{n+1}} 2^{-\ell}\ell^{-\varepsilon} \psi\big(2^{\ell}(x-2^{-(j-\ell_0)})\big)\\
 & & +\sum_{n=0}^{\infty} 2^{-j_{n+1}\alpha} \sum_{j=j_{n,\alpha}+1}^{j_{n+1}-1} \sum_{\ell=j_{n+1}}^{j_{n+1,\alpha}}\ell^{-\varepsilon}\psi\big(2^\ell (x-2^{-(j-\ell_0)})\big)
\end{eqnarray*}
satisfies the following properties:
\begin{enumerate}
\item $f$ is not a uniformly H\"olderian function,
\item the wavelet coefficients of $f$ satisfy property~(\ref{eq:wav-carWR}),
\item $f$ is uniformly irregular with exponent $\beta$, where
\begin{equation}\label{eq:beta}
 \beta=\max(\alpha\varepsilon,\frac{\alpha\varepsilon}{(1-\alpha)+\alpha\varepsilon})<\alpha.
\end{equation}
\end{enumerate}
\EPro \BProof The two first properties being straightforward, we just have to prove that $f$ is uniformly irregular with exponent $\beta$.
Let $n\in\N$ and define
\[
 f_j(x)=\sum_{\ell=j+2}^{j_{n,\alpha}} \ell^{-\varepsilon} \psi\big(2^\ell (x-2^{-(j-\ell_0)})\big),\\
\]
for $j\in\{j_n,\ldots,j_{n,\alpha}\}$ and
\[
 f_j(x)=\sum_{\ell=j+2}^{j_{n+1}} 2^{-\ell} \ell^{-\varepsilon} \psi\big( 2^\ell (x-2^{-(j-\ell_0)})\big)+
 \sum_{\ell=j_{n+1}}^{j_{n+1,\alpha}}\ell^{-\varepsilon}\psi\big(2^\ell (x-2^{-(j-\ell_0)})\big),
\]
for $j\in\{j_{n,\alpha},\ldots,j_{n+1}-1\}$. We need to estimate
\[
 f(2^{-(j-\ell_0)})-f(0) = f(2^{-(j-\ell_0)})
\]
for any $j\in\N$.
First, observe that for $j\neq j'$,  $\mathrm{supp}(f_j)\cap \mathrm{supp}(f_{j'})=\emptyset$. Indeed for any $j$, we have
\[
\mathrm{supp}(f_j)\subset [3.2^{-(j+2-\ell_0)},5.2^{-(j+2-\ell_0)}]\;
\]
and hence $f(2^{-(j-\ell_{0})})-f(0)=f_{j}(2^{-(j-\ell_{0})})$ for any $j\in\N$.

We now distinguish two cases.
Let us first assume that $j\in\{j_n,\ldots,j_{n,\alpha}\}$; we have
\[
 f(2^{-(j-\ell_0)})=
 2^{-j_n \alpha} \sum_{\ell=j+2}^{j_{n,\alpha}} \ell^{-\varepsilon} \psi(0)
 \ge 2^{-j_n \alpha}((j_{n,\alpha}+1)^{1-\varepsilon}-(j+2)^{1-\varepsilon})
\]
Therefore, if $j_n \le j\le j_{n,\alpha}/2$,
\[
 f(2^{-(j-\ell_{0})})
 \ge 2^{-j_n \alpha}(j_{n,\alpha}+1)^{1-\varepsilon}(1-2^{-(1-\varepsilon)})
 \ge C' 2^{-j\alpha\varepsilon},
\]
whereas if $j_{n,\alpha}/2\leq j\leq j_{n,\alpha}$,
\[
 f(2^{-(j-\ell_{0})})
 \ge 2^{-j_n \alpha}j_{n,\alpha}^{-\varepsilon}
 \ge j^{-1-\varepsilon}.
\]
Gathering these inequalities, we have, for any $j\in\{j_n, \ldots,j_{n,\alpha}\}$,
\begin{equation}\label{eq:ineqcex1a}
 f(2^{-(j-\ell_0)})\ge C' 2^{-j\alpha\varepsilon}.
\end{equation}

Let us now consider the second case, where $j\in\{j_{n,\alpha}+1,\cdots,j_{n+1}-1\}$ for some $n\in\N$. We have
\[
 f(2^{-(j-\ell_0)})= \big(2^{j_{n+1}(1-\alpha)}\sum_{\ell=j+2}^{j_{n+1}}2^{-\ell}\ell^{-\varepsilon}+ 2^{-j_{n+1}\alpha}\sum_{\ell=j_{n+1}}^{j_{n+1,\alpha}}\ell^{-\varepsilon}\big)\psi(0).
\]
If one remarks that
\begin{eqnarray*}
 f(2^{-(j-\ell_0)})
 &\ge& C' (2^{j_{n+1}(1-\alpha)} 2^{-j} j^{-\varepsilon} +2^{-j_{n+1}\alpha} j_{n+1,\alpha}^{1-\varepsilon}) \\
 &=& C' (2^{j_{n+1}(1-\alpha)}2^{-j} j^{-\varepsilon}+2^{-j_{n+1}\alpha\varepsilon}),
\end{eqnarray*}
then for any $j_{n,\alpha}+1\le j\le ((1-\alpha)+\alpha\varepsilon)j_{n+1}$, we get
\begin{equation}\label{eq:ineqcex1b}
 f(2^{-(j-\ell_{0})})
 \ge C' 2^{j\frac{1-\alpha}{(1-\alpha)+\alpha\varepsilon}} 2^{-j}j^{-\varepsilon}
 = C' 2^{-j\frac{\alpha\varepsilon}{(1-\alpha)+\alpha\varepsilon}}j^{-\varepsilon},
\end{equation}
whereas if $((1-\alpha)+\alpha\varepsilon)j_{n+1}\le j\le j_{n+1}-1$,
\begin{equation}\label{eq:ineqcex1c}
 f(2^{-(j-\ell_0)})
 \ge C' 2^{-j_{n+1}\alpha\varepsilon}
 \ge C' 2^{-j\frac{\alpha\varepsilon}{(1-\alpha)+\alpha\varepsilon}}.
\end{equation}
Inequalities~(\ref{eq:ineqcex1a}), (\ref{eq:ineqcex1b}) and (\ref{eq:ineqcex1c}) together imply $f\in UI^{\beta}(\R^{d})$.
\EProof

\subsection{Necessity of the logarithmic correction in the wavelet criteria}
Let $\varepsilon,\alpha\in (0,1)$, $\beta>1$ and define $(j_n)_{n\in\N}$ as
\[
 j_n=\floor{\beta^n},
\]
for any $n\in\N$. Let us also define the function $f_{\alpha,\beta,\varepsilon}$ on $\R$ as follows,
\begin{equation}\label{eq:fabe}
 f_{\alpha,\beta,\varepsilon}(x)=
 \sum_{n=0}^{\infty} \sum_{j=j_n+1}^{j_{n+1}} \frac{\inf(2^{-j_{n}\alpha},2^{j_{n+1}(1-\alpha)}2^{-j})}{j^\varepsilon} \sin(2^j \pi x).
\end{equation}
We first give an estimation of the wavelet coefficients $(c_{j,k})$ of $f_{\alpha,\beta,\varepsilon}$.
\BPro\label{pro:cex2a}
Assume that the multiresolution analysis is the Meyer multiresolution analysis. Then for $n\ge 1$, any $j\in\{j_n,\cdots,j_{n+1}-1\}$ and any $C>0$,
\begin{equation}\label{eq:wavcex2}
 \sup_{k\in\Z} |c_{j,k}|\le C\inf(2^{-j_n\alpha},2^{j_{n+1}(1-\alpha)}2^{-j}),
\end{equation}
for $n$ sufficiently large.
\EPro
\BProof
Let  $n\in\N$ and $\ell\in\{j_n,\ldots, j_{n+1}-1\}$. By definition of the wavelet coefficients of a bounded function, we have
\[
 c_{\ell,k}= 2^{\ell}\int_{\R^d} f_{\alpha,\beta,\varepsilon}(x) \psi(2^{\ell}x-k) \,dx.
\]
Since the trigonometric series $f_{\alpha,\beta,\varepsilon}$ is uniformly converging on any compact,
\[
 c_{\ell,k}=2^{\ell} \sum_{n=0}^{\infty} \sum_{j=j_n +1}^{j_{n+1}}\frac{\inf(2^{-j_n \alpha},2^{j_{n+1}(1-\alpha)}2^{-j})}{j^\varepsilon} \int_{\R^d} \sin(2^j \pi x)\psi(2^\ell x-k) \,dx,
\]
 or
\begin{eqnarray*}
 \lefteqn{c_{\ell,k}=} \\
 && \frac{2^\ell} {2i}\sum_{n=0}^{\infty}\sum_{j=j_n+1}^{j_{n+1}} \frac{\inf(2^{-j_n\alpha},2^{j_{n+1}(1-\alpha)}2^{-j})}{j^\varepsilon} \int_{\R^d} (e^{i 2^j \pi x}- e^{-i 2^j \pi x})\psi(2^\ell x-k)\,dx,
\end{eqnarray*}
that is,
\begin{eqnarray}\label{eq:sumcex2}
 \lefteqn{c_{\ell,k}=} \\
 && \sum_{n=0}^{\infty} \sum_{j=j_n+1}^{j_{n+1}} \frac{\inf(2^{-j_n \alpha},2^{j_{n+1}(1-\alpha)}2^{-j})}{j^\varepsilon} \frac{e^{i 2^{j-\ell}k\pi}\hat{\psi} (2^{j-\ell}k)-e^{-i2^{j-\ell}k\pi}\hat{\psi}(-2^{j-\ell}k)}{2i}. \nonumber
\end{eqnarray}
Since the Meyer wavelet belongs to the Schwartz class, its Fourier transform is symmetric and compactly supported with
\[
\mathrm{supp}(\hat{\psi})\subset [-\frac{8\pi}{3},-\frac{2\pi}{3}]\cup [\frac{2\pi}{3},\frac{8\pi}{3}],
\]
the sum in equality~(\ref{eq:sumcex2}) contains at most five terms corresponding to
\[
 k\in\{\ell-\log_2 (k),\ell-\log_2 (k)+1,\ell-\log_2(k)+2,\ell-\log_2(k)+3,\ell-\log_2(k)+4\}\;.
\]
One directly checks that for any $n\in\N$, $j\in\{j_n,\ldots, j_{n+1}-1\}$, this implies inequality~(\ref{eq:wavcex2}).
\EProof

Let us now prove the uniform irregularity properties of the functions $f_{\alpha,\beta,\varepsilon}$.
\BPro\label{pro:cex2b}
For any $\beta>1$ and any $(\alpha,\varepsilon)\in (0,1)^2$, $f_{\alpha,\beta,\varepsilon}\in UI^{\alpha}_{1-\varepsilon}(\R)$.
\EPro
\BProof
Let us remark that it is sufficient to prove that for any $\ell\in\N$,
\begin{equation}\label{eq:irregcex2}
 f_{\alpha,\beta,\varepsilon} (2^{-\ell})\ge 2^{-\alpha \ell} \ell^{1-\varepsilon}.
\end{equation}
Let $n_0\in\N$ and $\ell\in\{j_{n_{0}+1},\ldots,j_{n_0+1}\}$. By definition, we have
\begin{eqnarray*}
 f_{\alpha,\beta,\varepsilon} (2^{-\ell}) &=& \sum_{n=0}^{n_0-1}\sum_{j=j_n+1}^{j_{n+1}} \frac{\inf(2^{-j_n\alpha},2^{j_{n+1}(1-\alpha)}2^{-j})}{j^\varepsilon} \sin(2^j 2^{-\ell} \pi) \\
 && + \sum_{j=j_{n_0}+1}^{\ell-1} \frac{\inf(2^{-j_{n_0}\alpha},2^{j_{n_0+1}(1-\alpha)}2^{-j})}{j^\varepsilon}\sin(2^j\pi 2^{-\ell}).
\end{eqnarray*}
The classical inequality $\sin(x)\geq (2/\pi) x$ valid for any $x\in [0,\pi/2]$ leads to the following inequality if $j_{n_0}+1\le \ell\le j_{n_0+1}$,
\begin{eqnarray*}
 f_{\alpha,\beta,\varepsilon}(2^{-\ell}) &\ge& \sum_{j=j_{n_0}+1}^{\ell-1}\frac{\inf(2^{-j_n \alpha},2^{j_{n+1}(1-\alpha)}2^{-j})}{j^\varepsilon} 2^{j-\ell}\\
 &\geq& 2.2^{-\ell} \inf(2^{-j_{n_0}\alpha}2^\ell \ell^{1-\varepsilon},\ell^{1-\varepsilon} 2^{j_{n_0+1}(1-\alpha)})\\
 &\geq& 2\inf(\ell^{1-\varepsilon} 2^{-j_{n_0}\alpha},\ell^{1-\varepsilon}2^{-\ell} 2^{j_{n_0+1}(1-\alpha)}).
\end{eqnarray*}
Let $t\in (1,\beta)$ such that $\ell=tj_{n_{0}}$, that is $j_{n_{0}}=\ell/t$. We get
\[
 f_{\alpha,\beta,\varepsilon}(2^{-\ell})\ge
 2\inf(\ell^{1-\varepsilon}2^{-\ell\frac{\alpha}{t}},\ell^{1-\varepsilon} 2^{-\ell(1-\frac{\beta\ell}{t}+\frac{\alpha\beta\ell}{t})}).
\]
Since
\[
 \sup_{t\in [1,\beta]} \max(\alpha/t,1-\beta\ell/t+\alpha\beta\ell/t)\le \alpha,
\]
inequality~(\ref{eq:irregcex2}) is satisfied for any $\ell\in\N$.
\EProof

Propositions~\ref{pro:cex2a} and~\ref{pro:cex2b} together imply the following proposition.
\BPro\label{pro:cex2}
For any $(\alpha,\varepsilon,\beta)\in (0,1)^2\times (1,+\infty)$, the functions $f_{\alpha,\beta,\varepsilon}$ defined by the relation~(\ref{eq:fabe}) are uniformly H\"olderian, satisfy (\ref{eq:wav-carWR}) and belong to $UI^{\alpha}_{1-\varepsilon}(\R)$.
\EPro

\bibliographystyle{model1a-num-names}

\end{document}